\documentclass{scrartcl}

\KOMAoptions{paper=a4}
\KOMAoptions{fontsize=12pt}
\KOMAoptions{DIV=calc} 

\usepackage[utf8]{inputenc} 
\usepackage[T1]{fontenc} 
\usepackage[ngerman,english]{babel} 
\usepackage{csquotes} 
\usepackage{hyphenat} 

\usepackage{kpfonts}

\usepackage{cabin}

\usepackage{microtype} 

\addtokomafont{title}{\rmfamily}

\usepackage{enumitem}
\setlist[enumerate]{label*=(\alph*),ref=(\alph*),itemsep=0pt,topsep=5pt}
\setlist[itemize]{itemsep=0pt,topsep=5pt}

\KOMAoptions{DIV=calc}

\usepackage{booktabs} 

\usepackage{graphicx} 
\graphicspath{{pic/}} 

\usepackage[format=plain,labelfont={sf,footnotesize},textfont=small,margin=12pt]{caption} 
\DeclareCaptionLabelSeparator{MySpace}{\enskip } 
\newlength{\captionlength}

\usepackage[style=alphabetic,
						sorting=nyt,
						maxnames=4,
						backend=biber,
						doi=true,
						isbn=false,
						clearlang=true]
						{biblatex}

\bibliography{C:/Users/j.rau/Dropbox/Bibliography/biblio}

\DeclareFieldFormat*{title}{\textit{#1}} 
\renewbibmacro{in:}{} 
\DeclareFieldFormat{journaltitle}{\iffieldequalstr{journaltitle}{ArXiv e-prints}{Preprint}{#1}} 
\DeclareFieldFormat{booktitle}{#1} 
\DeclareRedundantLanguages{english,English}{english,german,ngerman,french}
\DeclareSourcemap{ 
  \maps[datatype=bibtex]{
    \map{
      \step[fieldset=addendum,null]
      \step[fieldset=eprintclass,null]
    }  
  }
}

\usepackage{hyperref} 
\usepackage{xcolor} 
\hypersetup{breaklinks=true} 
\definecolor{darkblue}{RGB}{0,0,170}
\definecolor{darkred}{RGB}{200,0,0}
\hypersetup{citecolor=darkblue, urlcolor=darkblue, linkcolor=darkblue, colorlinks=true} 
\usepackage[all]{hypcap} 

\usepackage{aliascnt} 
\usepackage[amsmath,thmmarks]{ntheorem} 

\newtheorem {theorem}{Theorem}[section]

\newaliascnt{proposition}{theorem}
\newtheorem {proposition}[proposition]{Proposition}
\aliascntresetthe{proposition}

\newaliascnt{lemma}{theorem}
\newtheorem {lemma}[lemma]{Lemma}
\aliascntresetthe{lemma}

\newaliascnt{corollary}{theorem}

\aliascntresetthe{corollary}

\newaliascnt{conjecture}{theorem}
\newtheorem {conjecture}[conjecture]{Conjecture}
\aliascntresetthe{conjecture}

\theorembodyfont{\upshape}
											 
\newaliascnt{definition}{theorem}
\newtheorem {definition}[definition]{Definition}
\aliascntresetthe{definition}

\newaliascnt{example}{theorem}
\newtheorem {example}[example]{Example}
\aliascntresetthe{example}

\newaliascnt{exercise}{theorem}

\aliascntresetthe{exercise}

\newaliascnt{goal}{theorem}

\aliascntresetthe{goal}

\newaliascnt{construction}{theorem}

\aliascntresetthe{construction}

\theoremheaderfont{\itshape}
\theorembodyfont{\upshape}

\newaliascnt{remark}{theorem}
\newtheorem {remark}[remark]{Remark}
\aliascntresetthe{remark}

\newaliascnt{convention}{theorem}
\newtheorem {convention}[convention]{Convention}
\aliascntresetthe{convention}

\newaliascnt{notation}{theorem}

\aliascntresetthe{notation}

\theoremstyle {nonumberplain}
\theoremheaderfont{\itshape}
\theorembodyfont{\upshape}
\theoremsymbol{\ensuremath{_\blacksquare}}

\newtheorem {proof}{Proof}

\usepackage {amsmath} 

\newcommand{\C}{{\mathbf C}}

\newcommand{\F}{{\mathbf F}}

\newcommand{\N}{{\mathbf N}}
\renewcommand{\P}{{\mathbf P}}

\newcommand{\R}{{\mathbf R}}
\newcommand{\T}{{\mathbf T}}

\newcommand{\Z}{{\mathbf Z}}

\newcommand{\CP}{\mathbf{CP}}

\newcommand{\CC}{{\mathcal C}}
\newcommand{\FF}{{\mathcal F}}
\newcommand{\GG}{{\mathcal G}}

\newcommand{\ZZ}{{\mathcal Z}}
\newcommand{\BM}{{\mathrm{BM}}}
\newcommand{\cell}{{\mathrm{cell}}}

\newcommand{\CCC}{{\check{C}}}
\newcommand{\HHH}{{\check{H}}}

\newcommand{\mm}{{\mathfrak m}}

\DeclareMathOperator{\OS}{OS}
\DeclareMathOperator{\Tr}{Tr}
\DeclareMathOperator{\rk}{rk}

\DeclareMathOperator{\val}{val}
\DeclareMathOperator{\id}{id}

\DeclareMathOperator{\Fix}{Fix}
\DeclareMathOperator{\fix}{fix}
\DeclareMathOperator{\Mat}{Mat}
\DeclareMathOperator{\gap}{gap}

\DeclareMathOperator{\Ver}{Vert}
\DeclareMathOperator{\cyc}{cyc}
\DeclareMathOperator{\relint}{relint}
\DeclareMathOperator{\im}{im}

\begin {document}

\title {The tropical Poincaré-Hopf theorem}
\author {Johannes Rau}
\date{} 

\maketitle

\begin{abstract}
	\noindent
  We express the beta invariant of a loopless matroid as tropical self-intersection number 
	of the diagonal of its matroid fan (a \enquote{local} Poincaré\hyp{}Hopf theorem). 
	This provides another example of uncovering the \enquote{geometry} of matroids by expressing
	their invariants in terms of tropicalised geometric constructions. 
	We also prove a global Poincaré-Hopf theorem and 
	initiate the study of a more general tropical Lefschetz\hyp{}Hopf trace formula
	by proving the two special cases of tropical curves and tropical tori.
\end{abstract}

{\footnote{{\noindent This work is part of the FAPA project \enquote{Matroids in tropical geometry} 
of the author at Universidad de los Andes, Colombia (pending)}}}

\section{Introduction}

The  Euler characteristic $\chi(X)$ of a compact manifold $X$ is equal to  the self-intersection number of the 
diagonal $\Delta_X \subset X \times X$. In short, $\deg \Delta_X^2 = \chi(X)$.
This is a reformulation  of the Poincaré-Hopf theorem \cite{Poi-SurLesCourbes,Hop-VektorfelderNDimensionalen}
(in view of $\deg \Delta_X^2 = \deg T_X$). 
We prove the following tropical (or matroid-theoretic) version:

\begin{theorem}[Local tropical Poincaré-Hopf] \label{mainresult}
  Let $M$ be a loopfree matroid of rank $n+1$. 
	Denote by $\Sigma_M = \Sigma'_M/\R\mathbf{1}$ its (projective) matroid fan and
	by $\Delta$ the diagonal of $\Sigma_M$ in $\Sigma_M \times \Sigma_M$. 
	Then the self-intersection of $\Delta$ in $\Sigma_M \times \Sigma_M$ is given by
	\begin{equation}  
    \deg \Delta^2 = (-1)^n \beta(M).
  \end{equation}
\end{theorem}

Here, $\beta(M)$ denotes the \emph{beta invariant} of a matroid.  
It is the canonical replacement for the 
Euler characteristic since e.g.\ $(-1)^n \beta(M) = \chi(U)$ 
if $U$ is the complement of a hyperplane arrangement realizing $M$
(also \autoref{lem:specialcase}).
The product $\Delta^2$ refers to the intersection product 
for tropical subcycles of matroid fans constructed in
\cite{Sha-TropicalIntersectionProduct, FR-DiagonalTropicalMatroid}. 
As mentioned before, we hope that this formula provides another
interesting instance of uncovering the geometric side of matroids 
by tropicalising a well-known classical geometric result.
Note, however, that our formula treats a non-compact
setup and has no immediate 
classical counterpart (c.f.\ \autoref{rem:surprise}). 
So our formula could be a starting point for finding 
similar classical statements (when $M$ is realizable).

The Poincaré-Hopf theorem can be regarded as a special case of the Lefschetz-Hopf trace formula (or fixed-point theorem) 
\cite[Section VII.6]{Lef-IntersectionsTransformationsComplexes,Dol-LecturesAlgebraicTopology}
applied to the endomorphism $\psi = \id_X$. 
Thus, \autoref{mainresult} naturally poses the question whether a 
more general tropical trace formula holds. 
Again, it is already interesting that such a statement can be formulated in the 
tropical setup without any compactness requirements, as follows (for details, see \autoref{sec:traceformula}).
Let $X$ be a smooth tropical variety of dimension $n$ without 
points of higher sedentarity (i.e.\ $X$ is locally isomorphic to (an open subset of) a matroid fan) 
and let $\psi \colon X \to X$ be a tropical endomorphism. 
We denote by $\Gamma_\psi$ and $\Delta$ the graph and diagonal, respectively, of $X$ in $X \times X$. 
The intersection product $\Gamma_\psi \cdot \Delta$ (again, in the sense of \cite{Sha-TropicalIntersectionProduct, FR-DiagonalTropicalMatroid})
can be regarded as  the cycle of \emph{stable} fixed points of $\psi$.
On the trace side, we use the tropical Hodge type homology groups $H_{p,q}(X)$ \cite{IKMZ-TropicalHomology, MZ-TropicalEigenwaveIntermediate}
and the Borel-Moore versions $H^{\BM}_{p,q}(X)$ (e.g.\ \cite{JRS-Lefschetz11Theorem}), both with real coefficients.
If $\psi$ is proper, we have induced pushforward maps $\psi_* \colon H^{\BM}_{p,q}(X) \to H^{\BM}_{p,q}(X)$.
We denote the traces by $\Tr(\psi_*, H^{\BM}_{p,q}(X))$. 

\begin{conjecture}[Tropical Lefschetz-Hopf trace formula] \label{traceformula}
  Let $\psi \colon X \to X$ be a proper tropical endomorphism of a smooth tropical variety $X$. Then we have
	\begin{equation} \label{eq:traceformula} 
		\deg (\Gamma_\psi \cdot \Delta) = \sum_{p,q} (-1)^{p+q} \Tr(\psi_*, H^{\BM}_{p,q}(X)).
	\end{equation}
\end{conjecture}

\autoref{mainresult} is a special case of this formula for $X = \Sigma_M$ and $\psi = \id$
(cf.\ \autoref{lem:specialcase}). 
In this paper, we also prove the following special cases.

\begin{theorem}[Global tropical Poincaré-Hopf] \label{GlobalPoincareHopf}
  The tropical Euler characteristic $\chi(X) := \sum_{p,q} (-1)^{p+q} \dim H_{p,q}(X)$ 
	of a smooth tropical variety without points of higher sedentarity
	is given by 
	\[
	  \chi(X) = \deg \Delta^2.
	\]
\end{theorem}

\begin{theorem}[Tropical Weil trace formula] \label{introcurves}
  If $X$ is a smooth tropical curve, \autoref{traceformula} holds.
\end{theorem}

\begin{theorem} \label{introtori}
  If $X$ is a tropical torus, \autoref{traceformula} holds.
\end{theorem}

In a work in progress we hope to prove \autoref{traceformula} in the case of 
\emph{matroidal} automorphisms (automorphisms which are induced by matroid automorphisms). 

\begin{remark} 
  Currently, \autoref{traceformula} should be restricted to varieties
	$X$ \emph{without points of higher sedentarity} since 
	we are lacking a definition for the intersection theoretic side in the presence
	of such points. 
	We hope that the intersection product (at least, its degree) can be defined 
	such that the statement holds in the greater generality. 
	In fact, in the case of curves a definition
	for points of higher sedentarity exists and 
	we will prove \autoref{introcurves} allowing such points.
	The tropical tori of \autoref{introtori} 
	do not contain such points by definition. 
\end{remark}

\begin{remark} \label{rem:surprise}
  As mentioned before, it might come as a surprise that the presented results 
	hold in the non-compact (even local) setup of e.g.\ \autoref{mainresult}. 
	Along these lines, we can make the following observation. 
  There is a canonically defined cycle class map $\cyc \colon Z_p(X) \to H_{p,p}^{\BM}(X)$ 
	which associates to a tropical $p$-dimensional subcycle its fundamental 
	cycle class. Moreover, the tropical homology groups for smooth varieties satisfy Poincaré duality and
	carry various intersection products 
	\cite{JSS-SuperformsTropicalCohomology, JRS-Lefschetz11Theorem, MZ-TropicalEigenwaveIntermediate}. 
	It is therefore tempting to hope that all the aforementioned statements can be proven
	using these constructions similar to proofs of the classical statement 
	(e.g.\ by writing down the Künneth decomposition for $\cyc(\Delta)$ 
	and showing that the intersection pairings are compatible). 
	This, however, does not work since	the tropical homology groups tend to be 
	too small in the non-compact setting.
	In particular, for matroid fans we have $H^{\BM}_{p,p}(\Sigma_M) = 0$ for $p \neq n$
	and thus $\cyc(\Delta)$ does not carry any information at all. 
	Hence, \autoref{mainresult} cannot be stated\fshyp{}proven 
	using intersection products on tropical homology alone. 
	In fact, our proof strategy has no classical analogue. 
	
	In the first version of this paper, I expressed my hope to 
	find classical counterparts of e.g.\ \autoref{mainresult}
	and to establish connections to the intersection theory\fshyp{}$K$-theory
	of (wonderful compactifications of) hyperplanes arrangements
	and to similar expressions for the characteristic polynomial of a matroid 
	in e.g.\ \cite{FS-kClassesMatroids, Alu-GrothendieckClassesChern, AHK-HodgeTheoryCombinatorial, LRS-ChernSchwartzMacpherson, ADH-LagrangianGeometryMatroids}.
	Since then, using their beautiful study of tautological classes of matroids, 
	a classical discussion of the realizable case has been given
	by Berget, Eur, Spink, Tseng in \cite[Appendix II]{BEST-TautologicalClassesMatroids}. 
\end{remark}

\begin{remark} \label{rem:otherhomology}
  A priori, \autoref{traceformula} could be formulated using any of the tropical (co)homology
	versions $H_{p,q}(X)$, $H^{\BM}_{p,q}(X)$, $H^{p,q}(X)$ and $H^{p,q}_c(X)$ 
	(for $H^{\BM}_{p,q}(X)$, we allow 
	locally finite chains; for $H^{p,q}_c(X)$, we restrict to compactly supported cochains). 
	We consider all these groups with real coefficients and hence drop $\R$ from the notation.
	By ordinary (not Poincaré) duality, we have $H_{p,q}(X) \cong (H^{p,q}(X))^*$ 
	as well as $H^{\BM}_{p,q}(X) \cong (H^{p,q}_c(X))^*$ and moreover $\psi_* = (\psi^*)^\top$
	(we assume $\psi$ proper in the second version). It follows that
	\begin{align} 
		\sum_{p,q} (-1)^{p+q} \Tr(\psi_*, H_{p,q}(X)) &= \sum_{p,q} (-1)^{p+q} \Tr(\psi^*, H^{p,q}(X)), \label{trace1} \\
		\sum_{p,q} (-1)^{p+q} \Tr(\psi_*, H^{\BM}_{p,q}(X)) &= \sum_{p,q} (-1)^{p+q} \Tr(\psi^*, H^{p,q}_c(X)).  \label{trace2}
	\end{align}
	So we are left with two possibilities for the trace side.
	If $\psi = \id$ (\autoref{GlobalPoincareHopf}), it follows from 
	Poincaré duality \cite{JSS-SuperformsTropicalCohomology, JRS-Lefschetz11Theorem}
	that all four versions do agree. Obviously, they also agree if $X$ is compact 
	as in \autoref{introtori}.
	For general endomorphisms, however, \autoref{trace1} and \autoref{trace2} 
	may be different and the 	$H_{p,q}(X)$ version may not give the correct answer. 
	For example, consider the
	standard tropical line $L\subset \R^2$ and the map $\psi \colon x \mapsto d x$, $d\in \N$. 
	Then the sum for $H_{*,*}(X)$ is $1-2d$, whereas for $H^{\BM}_{*,*}(X)$
	we get $d-2$. The latter number agrees with the intersection-theoretic side (see \autoref{curvelocalcase1}).
	This is why we decided to use $H^{\BM}_{p,q}$ in \autoref{traceformula}.
	On the other hand, in view \autoref{endomprop} and \autoref{workswithusualhom}
	one might argue that the given example is artificial and we could restrict 
	to automorphisms and use  $H_{p,q}$ without much loss.
	So, the question which version is more general\fshyp{}useful is probably still open for debate.
\end{remark}

\subsection*{Acknowledgements}

This project started during my visit to Oslo University in March 2019. 
My special thanks go to Kristin Shaw for the invitation, for bringing up this problem,
and for many useful discussions later on.
I would also like to thank the members of the tropical and matroidal seminar during semester
2020-1 at Universidad de los Andes for helpful discussions on the topic.
Finally, I would like to thank the reviewers for their carful reading and many helpful comments 
and corrections.

\section{The diagonal as complete intersection}

In this section, we mainly present a variant of a construction from \cite{FR-DiagonalTropicalMatroid}. It 
provides a description of the diagonal of a matroid fan in terms of $n$ tropical rational functions.
We will use this description later to compute the intersection product $\Delta^2$. 

\subsection{Preliminaries}

We start by fixing our basic notation for matroid fans. 
Throughout the paper, $M$ will denote a loopless matroid of rank $n+1$ on the ground set
$E=\{0, \dots, N\}$. 
We denote its rank function by $\rk$, its lattice of flats by $L(M)$ and the Möbius function
thereon by $\mu$. 
We recall that being \emph{loopfree} can be expressed as $\rk(\{i\}) = 1$ for all singletons $\{i\}$.

\begin{definition} 
  The \emph{beta invariant} of $M$ is 
	\[
	  \beta(M) := (-1)^{n+1} \sum_{F\in L(M)} \mu(\emptyset, F) \rk(F).
	\]
\end{definition}

We refer to \cite[Section 7.3]{Whi-CombinatorialGeometries} for more background on 
the beta invariant. It can be computed
asymmetrically as follows.

\begin{lemma} \label{BetaAsym}
  Fixing $0 \in E$, the beta invariant of $M$ is equal to 
  \[
	  \beta(M) = (-1)^n \sum_{\substack{F\in L(M) \\ 0 \notin F}} \mu(\emptyset, F).
	\]
\end{lemma}

\begin{proof}
  E.g.\ \cite[Proposition 7.3.1 (d)]{Whi-CombinatorialGeometries}
\end{proof}

Using $\beta(M)$ as \enquote{Euler characteristic} is motivated by the following well-known fact.  

\begin{proposition} 
  If $U$ is the complement of a complex hyperplane arrangement whose 
	associated matroid is $M$, then $\chi(U) = (-1)^n \beta(M)$.
\end{proposition}

\begin{proof}
	By the inclusion\fshyp{}exclusion properties of the Euler characteristic, we get
	\[
	  \chi(U) = \sum_{F\in L(M)} \mu(\emptyset, F) \chi(\C\P^{n - \rk(F)}) = \sum_{F\in L(M)} \mu(\emptyset, F) (n+1-\rk(F)).
	\]
	Since $\sum_{F\in L(M)} \mu(\emptyset, F) = 0$, the statement follows.
\end{proof}

To a loopless matroid $M$, we can associate an \emph{affine} matroid
fan $\Sigma'_M \subset \R^{N+1}$ whose lineality space contains the line $\R \mathbf{1}$,
see \cite[page 5]{Spe-TropicalLinearSpaces}.
Here, $\mathbf{1}$ denotes the all one vector $(1,\dots,1)$. 
By taking the quotient, we obtain the \emph{projective} matroid fan 
$\Sigma_M \subset \R^{N+1}/\R \mathbf{1} \cong \R^N$.
We identify $\R^N$ with $\R^{N+1}/\R \mathbf{1}$ by fixing the section $x_0 = 0$. 
We will mostly use the fan structure on $\Sigma_M$ associated to the lattice of flats $L(M)$,
called the \emph{fine subdivision} of $\Sigma_M$ \cite[Section 3]{AK-BergmanComplexMatroid}. 
This subdivision contains a cone $\sigma_\FF$ for any chain of flats $\FF$ as follows:
We will use the convention to write a chain of flats as a \emph{decreasing} sequence
\[
  E = F_0 \supsetneq F_1 \supsetneq \dots \supsetneq F_{l} \supsetneq F_{l+1} = \emptyset.
\]
We call $l(\FF) := l$ the \emph{length} of $\FF$. To such a chain, we associate the cone 
\[
  \sigma_\FF := \R_\geq \langle v_{F_1}, \dots, v_{F_l} \rangle + \R \mathbf{1} \subset \R^{N+1}.
\]
Here, for any subset $S \subset E$, we denote by $v_S \in \R^{N+1}$ the indicator vector for $S$
whose $i$-th entry is $-1$ if $i \in S$ (max-convention!) and $0$ if $i \notin S$. 
If no confusion is likely, we use the same notations $v_S$ and $\sigma_\FF$ for the projections to $\R^N$. 
Then $\Sigma'_M$ and $\Sigma_M$ can be described as to the union of all such cones
in $\R^{N+1}$ and $\R^N$, respectively. 
Note that $\Sigma_M$ is a unimodular fan of pure dimension $n$. 
We will use the same notation $\Sigma_M$ for the collection of cones as well as the underlying set.

Using chains $\CC$ of arbitrary subsets $C_i \subset E$, the cones $\sigma_\CC$ form a subdivision
of $\R^N$ (and $\R^{N+1}$) called the \emph{braid arrangement fan}. 
It can be equivalently described as the intersection of the hyperplane 
subdivisions $x_i = x_j$ for all $i \neq j$ or as the normal fan
of the permutahedron. The fine subdivisions of matroid fans are subfans of the
braid arrangement fan. 

Following \cite{Spe-TropicalLinearSpaces}, we recall that to any point $x \in \R^N$
we can associate a matroid $M_x$ whose bases are the $x$-maximal bases of $M$.
Here, the $x$-weight of a basis $B$ is defined to be $\langle x, v_B \rangle$.
A point $x$ is contained in $\Sigma_M$ if and only if $M_x$ is loopfree. 
Using this description, it is easy to show that 
$\Sigma'_{M_1 \oplus M_2} = \Sigma'_{M_1} \times \Sigma'_{M_2}$.
We are interested in the variant of this statement for projective fans. 

The analogous result for projective matroid fans in based 
on the notion of parallel connection of $M_1$ and $M_2$ which we quickly recall.
Let $E_1$ and $E_2$ be the ground sets of $M_1$ and $M_2$, respectively, and assume
$0 \in E_1 \cap E_2$. The \emph{parallel connection} of $M_1$ and $M_2$ along $0$ 
is the matroid $M = M_1 \oplus_0 M_2$ whose ground set $E$ 
is the disjoint union of $E_1$ and $E_2 \setminus \{0\}$. 

\begin{convention}
  Throughout the following, 
	we will describe subsets of $E$ in a more symmetric fashion as tuples $(S_1,S_2)$, 
	$S_1 \subset E_1$, $S_2 \subset E_2$ such that either $0 \in S_1 \cap S_2$
	or $0 \notin S_1 \cup S_2$ (the associated subset is $S_1 \sqcup (S_2 \setminus \{0\})$). 
\end{convention}

Using this convention, $M$ is defined as the matroid
whose flats are the pairs $(F_1, F_2)$ such that 
$F_i \in L(M_i)$ for $i=1,2$. 
It follows that the rank of such a flat 
is 
\[
  \rk_M(F_1, F_2) = \rk_{M_1}(F_1) + \rk_{M_2}(F_2) - \delta_{0 \in F_1}
\]
with $\delta_{0 \in F_1} = 1$ if $0 \in F_1$ and $\delta_{0 \in F_1} = 0$ otherwise.
In terms of bases, we have the following description: 
A basis of $M$ is a pair $(B_1, B_2)$
such that $0 \in B_1 \cap B_2$ and $B_1$ and $B_2$ are bases of
$M_1$ and $M_2$, respectively, or $0 \notin B_1 \cup B_2$
and $B_1$ and $B_2 \cup \{0\}$ are bases of
$M_1$ and $M_2$, or the symmetric version of the second case. 
We have a canonical identification of ambient spaces
$\R^{E\setminus 0} = \R^{E_1 \setminus 0} \times \R^{E_2 \setminus 0}$
which is compatible with our convention of setting $x_0 =0$.

\begin{lemma} \label{lem:cartproduct}
  Let $M$ be a parallel connection of two loopfree matroids $M_1$ and $M_2$. 
	Then $\Sigma_M = \Sigma_{M_1} \times \Sigma_{M_2}$ (as sets). 
\end{lemma}

\begin{proof}
  From the above description of bases it is clear that the 
	$(x,y)$-weight of  a basis $(B_1, B_2)$ of $M$ is equal to
  $x$-weight of $B_1$ plus the $y$-weight of $B_2$ 
	(again, recall that $x_0 = 0$ and $y_0=0$ by convention).
	It follows that 
	\begin{itemize}
		\item the element $0$ is contained in a basis of $M_{(x,y)}$ if and only if it 
		      is contained in a basis of $(M_1)_x$ and a basis of $(M_2)_y$,
	  \item in this case, $M_{(x,y)}$ is equal to the parallel connection 
		      of $(M_1)_x$ and $(M_2)_y$ along $0$.  
	\end{itemize}
	In particular, $M_{(x,y)}$ is loopfree if and only if both
	$(M_1)_x$ and $(M_2)_y$ are loopfree, which proves the claim. 
\end{proof}

\subsection{Generic chains of matroids}

Given two matroids $M, N$ on the ground set $E$, it is obvious that $\Sigma_N \subset \Sigma_M$
(both as sets and fans) if and only if $L(N) \subset L(M)$ (or, in matroid terminology, $N$ is a \emph{quotient} of $M$). 
In such a case, there exists a canonical sequence of matroids $N = M_0, M_1, \dots, M_s = M$
such that $\rk(M_i) = \rk(N) + i$ and 
\begin{equation} \label{eq:genericmatroids} 
  \Sigma_N = \Sigma_{M_0} \subset \Sigma_{M_1} \subset \dots \subset \Sigma_{M_s} = \Sigma_{M},
\end{equation}
see \cite[Corollary 3.6]{FR-DiagonalTropicalMatroid}. 
These matroids are given in terms of their rank functions by
\begin{equation} \label{eq:intermediatematroids} 
  \rk_{M_i}(S) = \min\{\rk_N(S) + i, \rk_M(S)\}.
\end{equation}
If $N$ and $M$ correspond to hyperplane arrangements associated
to the projective subspaces $K \subset L \subset \CP^n$,
then the $M_i$ correspond to a chain of generic subspaces
$K = S_0 \subset S_1 \subset \dots \subset S_s = L$.
Moreover, there exists an associated sequence of rational functions
$g'_1, \dots, g'_s : \R^{N+1} \to \R$ such that
\begin{equation} \label{eq:intersectfunctiongeneral} 
  \Sigma'_{M_{s-i}} = g'_i \cdot g'_{i-1} \cdots g'_1 \cdot \Sigma'_M,
\end{equation}
see \cite[Proposition 3.10]{FR-DiagonalTropicalMatroid}
(we refer to \cite[Definition 3.4]{AR-FirstStepsTropical} for a definition of the
intersection with rational functions\fshyp{}the divisor construction). 
These functions are linear on each cone of the braid arrangement fan of $\R^{N+1}$ and
hence determined by their values on the indicator vectors $v_S$, $S \subset E$. 
These values are given by
\begin{equation} \label{eq:functioningeneral} 
  g'_i(S) = \begin{cases}
	           -1 & \rk_M(S) \geq \rk_N(S) + s + 1 - i, \\
						  0 & \text{otherwise}.
	         \end{cases}
\end{equation}
Here, we use the short-hand $g'_i(S)$ instead of $g'_i(v_S)$. 
We are now going to use this construction in the special case of diagonals.

\subsection{Cutting out the diagonal}

The construction from \autoref{eq:intersectfunctiongeneral}
was used in \cite{FR-DiagonalTropicalMatroid} to construct
functions that cut out the diagonal of a matroid fan. This
description was then used to define a general intersection
product for tropical subcycles of matroid fans (and smooth tropical
varieties without points of higher sedentarity). 
We will use a slight variant of this construction here. 
Let us quickly explain the difference. 

In \cite{FR-DiagonalTropicalMatroid}, the construction of 
\autoref{eq:genericmatroids} was applied to the \emph{affine} matroid fans,
i.e.\ to the diagonal $\Delta'$ of $\Sigma'_M$ in $\Sigma'_M \times \Sigma'_M = \Sigma'_{M \oplus M}$. 
This is not quite what we want, since we are interested in the self-intersection of \emph{projective} fans
$\Delta \subset \Sigma_M \times \Sigma_M$
(the self-intersection of $\Delta' \subset \Sigma'_{M \oplus M}$ is always zero). We could consider the lift 
$\widetilde{\Delta} \subset \Sigma'_M \times \Sigma'_M$ of $\Delta$ and compute the self-intersection
$\widetilde{\Delta}^2$.
But note that $\widetilde{\Delta} = \Delta' + \R (\mathbf{1}, \mathbf{0})$ has \emph{extra} lineality
space. In particular, unlike $\Delta'$ and $\Delta$, it is not a matroid fan (for the given embedding) 
and it is not clear to the author how to express 
$\widetilde{\Delta}$ as complete intersection 
or how to compute the self-intersection $\widetilde{\Delta}^2$ otherwise. 
Instead, we will work directly with the projective matroid fans
(and justify why the results agree). 
This requires to break symmetry by choosing an element in $E$. However, we will later
see that this corresponds nicely to the asymmetric formula for $\beta(M)$ in
\autoref{BetaAsym}. 
Let us give the details here.

We denote by $M \oplus_0 M$ the parallel connection $M$ with it itself along $0$. 
%
%
By \autoref{lem:cartproduct}, 
\[
  \Sigma_M \times \Sigma_M = \Sigma_{M \oplus_0 M}.
\]

We can describe $\Delta \subset \Sigma_M \times \Sigma_M$ as the matroid fan given by the rank function
\begin{equation} 
  \rk_\Delta(F,G) = \rk(F \cup G).
\end{equation}
Of course, the union $F \cup G$ is understood to be the (non-disjoint) union in $E$ here. 
It is straightforward to check that this indeed defines a matroid $M_\Delta$ and 
that the flats of $M_\Delta$ are of the form $(F,F)$, 
$F \in L(M)$, which shows that $\Delta = \Sigma_{M_\Delta}$ 
(even on the level of fan structures).

We now apply the construction of \autoref{eq:genericmatroids} to 
$M_\Delta$ and $M \oplus_0 M$ (of rank $n+1$ and $2n+1$, respectively). 
From \autoref{eq:functioningeneral}
we get functions $g'_1, \dots, g'_n \colon \R^{2N+1} \to \R$ given by
\begin{equation} \label{eq:functiongi} 
  g'_i(F,G) = \begin{cases}
	             -1 & \rk(F) + \rk(G) - \delta_{0 \in F} \geq \rk(F \cup G) + n + 1 - i, \\
						    0 & \text{otherwise}.
					   \end{cases}
\end{equation}
Note that these functions are \enquote{homogeneous} functions in the sense that they
live on $\R^{2N+1}$ (not $\R^{2N}$) and $g'_i(\mathbf{1}) = - g'_i(E,E) = 1$. 
In order to obtain functions on $\R^{2N}$, we break the symmetry again 
and dehomogenise these functions by subtracting the coordinate function $x_0$. Then the functions
$g_i = g'_i - x_0 \colon \R^{2N} \to \R$ are well-defined and determined by the values
\begin{equation} \label{eq:functiongidehom} 
  g_i(F,G) = \begin{cases}
	             -1 & 0 \notin F, \rk(F) + \rk(G) \geq \rk(F \cup G) + n + 1 - i, \\
	             +1 & 0 \in F, \rk(F) + \rk(G) \leq \rk(F \cup G) + n + 1 - i, \\
						    0 & \text{otherwise}.
					   \end{cases}
\end{equation}
Our discussion so far can be summarized in the following statement, which follows
directly from \autoref{eq:intersectfunctiongeneral}.

\begin{proposition} \label{prop:diagonalintersection}
  The diagonal $\Delta$ of $\Sigma_M$ in $\Sigma_M \times \Sigma_M$ can be described
	as the complete intersection
	\[
	  \Delta = g_n \cdots g_1 \cdot (\Sigma_M \times \Sigma_M)
	\]
	using the rational functions $g_i : \R^{2N} \to \R$ from \autoref{eq:functiongidehom}.
\end{proposition}

Note that we use the term \enquote{complete intersection} in a rather weak sense here, since 
in general the functions $g_1, \dots, g_n$ are not tropically polynomial nor otherwise
convex\fshyp{}regular. 

In order to compute $\Delta^2$, we will now just restrict the functions to $\Delta$,
or rather, consider the pullbacks along $d \colon x \mapsto (x,x)$. 
Set $f_i := d^*(g_i) \colon \R^N \to \R$. Note that $d$ is compatible with the braid arrangement fans of $\R^N$
and $\R^{2N}$ (i.e.\ it maps a cone $\sigma_\CC$ to the cone $\sigma_{\CC'}$ where
$\CC'$ is the chain of subsets obtained from replacing each $C$ in $\CC$ by $(C,C)$).
Hence the functions $f_i$ are linear on the cones of the braid arrangement fan
on $\R^N$ and completely determined by the values 
\begin{equation} \label{eq:functionsfi} 
  f_i(F) = \begin{cases}
	             -1 & 0 \notin F \text{ and } \rk(F) \geq n + 1 - i, \\
							 +1 & 0 \in F \text{ and } \rk(F) \leq n + 1 - i, \\
						    0 & \text{otherwise}.
					   \end{cases}
\end{equation}
Summarizing again, we can describe the left hand side of \autoref{mainresult}
as follows.

\begin{proposition} \label{mainresult2}
  Let $\Delta \subset \Sigma_M \times \Sigma_M$ by the (projective) 
	diagonal of a matroid fan of loopless matroid $M$. 
	Then
  \begin{equation}  
    \deg \Delta^2  = \deg(f_n \cdots f_1 \cdot \Sigma_M).
  \end{equation}
\end{proposition}

\begin{proof}
  By \autoref{prop:diagonalintersection} we can express $\Delta$ as $g_n \cdots g_1 \cdot (\Sigma_M \times \Sigma_M)$.
	Hence $\deg \Delta^2 = \deg(g_n \cdots g_1 \cdot \Delta)$ by \cite[Theorem 4.5 (6)]{FR-DiagonalTropicalMatroid}.
	Finally, by the projection formula \cite[Proposition 7.7]{AR-FirstStepsTropical} it follows that the latter expression
	is equal to $\deg(f_n \cdots f_1 \cdot \Sigma_M)$.
\end{proof}

\section{The local Poincaré-Hopf theorem}

We will now give a description of the intermediate intersection products
$f_k \cdots f_1 \cdot \Sigma_M$ and prove the description by induction. 
\autoref{mainresult} will then just follow by inspection of the case $k = n$. 
The intermediate intersection products can be described explicitly, 
and it would be interesting to see if they appear in other situations or if they can be related
to other canonical tropical subcycles of $\Sigma_M$ such as the CSM classes 
from \cite{LRS-ChernSchwartzMacpherson}
or the intersection products appearing in 
\cite{AHK-HodgeTheoryCombinatorial,ADH-LagrangianGeometryMatroids}. 

We set $X_k := f_k \cdots f_1 \cdot \Sigma_M$. In order to describe the $X_k$, let 
us introduce some terminology first. 
We denote the rank function of $M$ by $\rk$. Let $\FF$ be a chain of flats 
$E = F_0 \supsetneq F_1 \supsetneq \dots \supsetneq F_{l} \supsetneq F_{l+1} = \emptyset$.
Note that since $X_k$ is of codimension $k$ in $\Sigma_M$, it is a union 
of cones $\sigma_\FF$ where $\FF$ is a chain obtained by removing $k$ flats from 
a maximal chain in $L(M)$. 

\begin{definition} \label{def:gapsequence}
  The \emph{(rank) gap sequence} of $\FF$, denoted by $\gap(\FF) = (r_0, \dots, r_l)$, 
	is the sequence of numbers $r_i := \rk F_{i} - \rk F_{i+1} -1$.
\end{definition}

As an example, consider the uniform matroid $U_{3,4}$. 
Then the gap sequence for $E \supsetneq \{0,1\} \supsetneq \emptyset$
is $(0,1)$ whereas the gap sequence for 
$E \supsetneq \{0\} \supsetneq \emptyset$ is $(0,1)$. 

To describe $X_k$, we will only need chains $\FF$ 
whose gap sequences have one of the following
two shapes:
\begin{equation} \begin{split} \label{eq:rankgaps} 
  \gap(\FF)  &= (r, s, 0, \dots, 0), \;\; r+s = k, \;\; r,s \geq 0, \\
	\gap(\FF)  &= (k, 0, \dots, 0).
\end{split} \end{equation} 
More specifically, we are only interested in the following two cases
(c.f.\ \autoref{chaintypes}). 

	\begin{figure}[!b]
		\centering
		\input{pic/chaintypes.TpX}
		\caption{Chains of type $(r,s)$ and type $(k)$}
		\label{chaintypes}
	\end{figure}%

\begin{definition} 
  A chain $\FF$ is of \emph{type $(r,s)$} if its gap sequence is 
	$(r, s, 0, \dots, 0)$ and additionally $0 \notin F_1$
	(i.e.\ the only term of $\FF$ containing $0$ is $E$).
	Here, we include the trivial chain $E \supsetneq \emptyset$
	which in our terminology is the unique chain of type $(n,0)$. 

  A chain $\FF$ is of \emph{type $(k)$} if its gap sequence is 
	$(k, 0, \dots, 0)$ and additionally $0 \in F_1$.
\end{definition}

Note that the types $(k,0)$ and $(k)$ have
the same gap sequence, but differ as to whether $0$ is contained in $F_1$ or not. 
By extension, the type of a cone in $\Sigma_M$ 
is the type of the corresponding chain (if of any type at all). 
We can now describe the cycles $X_k$ (c.f.\ \autoref{exampleX1}).

\begin{proposition} \label{prop:descriptionXk}
  The tropical cycle $X_k$ consists of the cones of type
	$(r,s)$, $r+s=k$, $r,s \geq 0$, and the cones of type $(k)$. 
	The weight of a cone $\sigma$
	of type $(r,s)$ in $X_k$ is given by 
	\begin{equation} \label{eq:weightsXk} 
		\omega(\sigma) = (-1)^{\rk(M/F_1) - 1} \beta(M/F_1). 
	\end{equation}
	The weight of a cone $\sigma$
	of type $(k)$ in $X_k$ is always $\omega(\sigma) = 1$.
\end{proposition}

Before proving the proposition, let us 
check consistency by showing that the case $k=n$ 
implies \autoref{mainresult}. 

\begin{proof}[\autoref{mainresult}]
  In the case $k=n$, the only chain of correct dimension is the 
	trivial chain $E \supset \emptyset$ of type $(n,0)$. 
	By \autoref{prop:descriptionXk}, its weight is $(-1)^n \beta(M)$.
\end{proof}

  \setlength{\captionlength}{0.9\textwidth}
  \addtolength{\captionlength}{-0.45\textwidth}
	\begin{figure}[b]
		\begin{minipage}[r]{\captionlength}
			\caption{The cycle $X_1$ in the case the standard hyperplane in $\R^3$ (i.e.\ the uniform matroid $U_{3,4}$
                                  of rank $3$ on $4$ elements). The rays of type $(1)$, $(1,0)$ and $(0,1)$ are displayed in black, blue and red,
																	respectively.}
			\label{exampleX1}
		\end{minipage} \hfill	
		\begin{minipage}[l]{0.45\textwidth}
			\input{pic/exampleX1.TpX}
		\end{minipage}%
	\end{figure}

\begin{remark} 
	The appearance of the beta invariant of the \enquote{factors}
	of the chain $\FF$ is reminiscent of
	the definition of tropical CSM cycles in \cite[Definition 5]{LRS-ChernSchwartzMacpherson}.
	However, CSM cycles take into account cones with arbitrary gap sequences
	and the beta invariant of each factor in the chain. 
	In contrast, in $X_k$ only special gap sequences and the beta invariant
	of the first factor occur. For $k=n$ the differences disappear and 
	\autoref{mainresult} can also be stated as 
	\[
	  \Delta^2 = \text{csm}_0(\Sigma_M).
	\]
\end{remark}

We now want to prove \autoref{prop:descriptionXk}
using, of course, an induction on $k$. 
In order to perform the step $k \to k+1$, 
let us first study the codimension one faces of $X_k$. 

\begin{lemma} \label{lem:balancingXk}
  Assume that \autoref{prop:descriptionXk} holds for $X_k$ and
	let $\GG$ be a chain of flats corresponding to a codimension one
	face $\tau$ of $X_k$.
	Let $G \supsetneq H$ denote the flats in $\GG$ 
	corresponding to the last non-zero entry of $\gap(\GG)$. Then the following holds:
	\begin{enumerate}
		\item 
		  Exactly one of the following four statements holds true. 
			\begin{align*} 
				\text{(A)} &&& \rk(G \cup \{0\}) \leq n-k & \text{(B)} &&& \rk(G) \geq n-k, 0 \notin G,  \\
				\text{(C)} &&& G = E, 0 \in H,            & \text{(D)} &&& G=E,  0 \notin H. 
			\end{align*}
		\item
		  Assume that $\GG$ belongs to one of the cases (A), (B) or (C). Then the facets of $X_k$ containing $\tau$
			correspond bijectively to fillings $G \supsetneq F \supsetneq H$ with flats
			$F$ of rank $\rk F = \rk H + 1$. Moreover, all such facets have identical type and weight $\omega$
			and $X_k$ is balanced at $\tau$. Explicitly, the balancing condition at $\tau$ is given by
			\begin{equation} \label{eq:balancingA} 
				\sum_F \omega v_F = \omega v_G + \omega (\val - 1) v_H,
			\end{equation}
			where $F$ runs through all such fillings and $\val$ denotes the number of fillings.
	\end{enumerate}
\end{lemma}

\begin{proof}
  Let $\FF$ be a chain corresponding to a facet of $X_k$ containing $\tau$. By definitions,
	$\GG$ is obtained from $\FF$ by removing one of its flats, say $F_i$. 
	If $i > 2$ or if $\FF$ is of type $(k)$ and $i > 1$, it follows that $\GG$ satisfies (A). 
	If $\FF$ is of type $(r,s)$ and $i = 2$, we obtain case (B). Finally, if $i = 1$,
	we end up with (C) or (D), depending on whether $0 \in F_2$ or not. It is clear that
	the cases are mutually exclusive. Hence (a) follows.
	
	For (b), we have a closer look at the previous argument. Note that the cases (A), (B), (C)
	correspond exactly to the case where either $\FF$ is of type $(r,s)$ and $i \geq 2$
	or $\FF$ is of type $(k)$ and $i \geq 1$. In both cases, $\rk(F_i) = \rk(F_{i+1}) + 1$
	and hence $\FF$ corresponds to a filling as described in the statement. Moreover, 
	it is obvious that each filling occurs in this way. The type and weight of $\FF$ is
	completely determined by the principal part $E \supsetneq F_1$ of $\FF$
	which is still present in $\GG$ and hence
	is fixed for given $\GG$. It remains to check the balancing condition in the form
	of \autoref{eq:balancingA}. 
	After dividing by $\omega$, this follows from the well-known fact
	that the sets $F \setminus H$, running through flats $F$ with 
	$G \supsetneq F \supsetneq H$ and $\rk F = \rk H + 1$, 
	form a disjoint partition of $G \setminus H$.
\end{proof}

We are now ready to prove \autoref{prop:descriptionXk}. 

\begin{proof}[\autoref{prop:descriptionXk}]
  The induction start $k=0$ is trivial (note that $\beta(M)=1$ for any 
	loopless matroid of rank $1$).
	Let us prove the step $k \to k+1$. 
	For each codimension one cone $\tau$ of $X_k$, we need to compute
	its weight in $f_{k+1} \cdot X_k$. 
	For convenience, we recall from \autoref{eq:functionsfi} that $f_{k+1}$ is given by
	\begin{equation}  
		f_{k+1}(F) = \begin{cases}
								 -1 & 0 \notin F \text{ and } \rk(F) \geq n - k, \\
								 +1 & 0 \in F \text{ and } \rk(F) \leq n - k, \\
									0 & \text{otherwise}.
							 \end{cases}
	\end{equation}
	Copying the notation from \autoref{lem:balancingXk}, we
	denote by $\GG$ the chain associated to $\tau$ and
	by $G \supsetneq H$ the last non-trivial step in $\GG$. 
	We go through the cases (A), (B), (C), (D) according
	to \autoref{lem:balancingXk} (a). 
	By \autoref{eq:balancingA}, the weight of $\tau$ in in $f_{k+1} \cdot X_k$ is 
	\begin{equation} \label{eq:weightStep} 
		\omega \cdot \left( \sum_F  f_{k+1}(F) -  f_{k+1}(G) - (\val - 1) f_{k+1}(H)\right).
	\end{equation}

	Let $r,s$ denote integers 
	such that $r+s=k$ and $r,s \geq 0$. 
	Note that $\GG$ is of type $(k+1,0)$ in case (D),
	of type $(k+1)$ in case (C), 
	of type $(r,s+1)$ in case (B),
	and of different shape (not present in the description 
	of $X_{k+1}$ for case (A). 
	
	If $\GG$ is of type (A), then all the flats $F$ occurring in 
	\autoref{eq:balancingA} (including $G$ and $H$) satisfy $\rk(F \cup \{0\}) \leq n-k$.
	Hence, their value under $f_{k+1}$ is $+1$ when $0 \in F$ and $0$ if not.
	Plotting these values, we get
	the following three possible shapes.
	\begin{align*} 
			                        & 0 \in H   && 0 \in G \setminus H && 0 \notin G \\ 
			f_{k+1}(G) \hspace{2em} & 1         && 1                   && 0 \\
			f_{k+1}(F) \hspace{2em} & 1 \dots 1 && 1 0 \dots 0         && 0 \dots 0 \\ 
			f_{k+1}(H) \hspace{2em} & 1         && 0                   && 0
	\end{align*}
	(Note that the single $1$ in the middle of the $f_{k+1}(F)$ line corresponds to 
	$\overline{H \cup 0}$). 
	By \autoref{eq:weightStep} the weight of $\tau$ in $f_{k+1} \cdot X_k$ is 
	zero in all three cases.

	Let us now assume $\GG$ satisfies (B). Then the ranks of the flats involved in 
	\autoref{eq:balancingA} are $\rk(G) \geq n-k$, $\rk(F) = n-k-1$ and $\rk(H) = n-k-2$,
	and none of these flats contains $0$. 
	Hence the pattern of values under $f_{k+1}$ is:
	\begin{align*} 
			f_{k+1}(G)  \hspace{2em} & 1         \\ 
			f_{k+1}(F)  \hspace{2em} & 0 \dots 0 \\ 
			f_{k+1}(H)  \hspace{2em} & 0         
	\end{align*}
	Hence the weight assigned to $\tau$ is $\omega = (-1)^{r} \beta(M/F_1)$, 
	as required.
	
	We continue with case (C), so now $E=G$ and $0 \in H = F_1$. The ranks are now given
	by $\rk(G) = n+1$, $\rk(F) = n-k$ and $\rk(H) = n-k-1$. Since all flats contain $0$,
	we get the pattern of values
	\begin{align*} 
		f_{k+1}(G)  \hspace{2em} & 0         \\ 
		f_{k+1}(F)  \hspace{2em} & 1 \dots 1 \\ 
		f_{k+1}(H)  \hspace{2em} & 1         
	\end{align*}
  which gives weight $1$ (note that $\omega = 1$ in this case).
	
	Finally, we are left with case (D). So now $E=G$, $0 \notin H = F_1$
	and $\gap(\GG) = (k+1, 0, \dots, 0)$. 
	It follows from the description of $X_k$ that 
	the facets of $X_k$ containing $\tau$ correspond
	to $F = \overline{H \cup 0}$ and to all flats 
	$F \supsetneq H$ with $0 \notin F$. 
	Note that the balancing condition around
	$\tau$ written in terms of the vectors primitive generators $v_F$
	may only involve the additional vectors $v_E = \mathbf{1}$ and $v_H$
	(with certain coefficients). 
	But $f_{k+1}(H) = f_{k+1}(E) = 0$ (since $\rk H = n-k-1$), which
	means we can compute the weight of $\tau$ without knowing the coefficients
	(in fact, it can be checked that they are both equal to $1$). 
	To do so, note that $\rk(\overline{H \cup 0}) = n - k$ and $\rk(F) \geq n-k$
	for all the flats $F \supsetneq H$ with $0 \notin F$, so 
	they all evaluate to $1$ under $f_{k+1}$. Moreover, the weight
	of the facet associated to $\overline{H \cup 0}$ is $1$, while
	for the other flats $F$ it is $(-1)^{\rk(M/F) - 1} \beta(M/F)$.
	Hence the weight of $\tau$ in $f_{k+1} \cdot X_k$ is
	equal to 
	\[
	  \omega(\tau)  = 1 - \sum_{\substack{F \supsetneq H \\ 0 \notin F}} (-1)^{\rk(M/F) - 1} \beta(M/F).
	\]
	To finish the calculation, we recall that $\mu$ can be defined as the inverse of the zeta function of $L(M)$
	and hence for any interval $H \subset G$ satisfies  
	\[
	  \sum_{\substack{F \in L(M) \\ H \subset F \subset G}} \mu(F, G) = \delta(H,G),
	\]
	where the delta fucntion is $\delta(H,G) = 0$ unless $H=G$, in which case $\delta(H,G) = 1$. 
	Using  \autoref{BetaAsym} twice, we can now compute $\omega(\tau)$ as follows.
	\begin{equation*} \begin{split} 
		\omega(\tau) &= 1 - \sum_{\substack{F \supsetneq H \\ 0 \notin F}} (-1)^{\rk(M/F) - 1} \beta(M/F) 
								 = 1 - \sum_{\substack{F \supsetneq H \\ 0 \notin F}} \sum_{0 \notin G} \mu(F,G) \\
								 &= 1 - \sum_{0 \notin G} \sum_{H \subsetneq F \subset G} \mu(F,G) 
								 = 1 - \sum_{0 \notin G} (\delta(H,G) - \mu(H,G)) \\
								 &= 1 - 1 + \sum_{0 \notin G}  \mu(H,G) 
								 = (-1)^{\rk(M/H) - 1} \beta(M/H).
	\end{split} \end{equation*}
  This agrees with \autoref{eq:weightsXk}, so we are done.
\end{proof}

\section{Towards a tropical Lefschetz-Hopf trace formula} \label{sec:traceformula}

In this section, we want to discuss three special cases of \autoref{traceformula}. 
A \emph{global} Poincaré-Hopf theorem corresponding to $\psi = \id$, the case of tropical curves
(even with points of higher sedentarity) and the case of tropical tori. 
The follow-up paper \cite{Rau-TropicalLefschetzHopf}
contains a proof of \autoref{traceformula} in the case of 
\emph{matroidal} automorphisms (automorphisms which are induced by matroid automorphisms). 
We start by giving details concerning the 
notions used in the introduction.

\subsection{Smooth tropical varieties}

Let $\Sigma$ be a polyhedral fan in $\R^N$. By abuse of notation, we use the same letter 
$\Sigma$ to the denote the support of $\Sigma$. 
For $x \in \Sigma$, we set 
\[
	U(x) := \bigcup_{\substack{\sigma \in \Sigma \\ x \in \sigma}} \relint(\sigma).
\]
A open subset $U \subset \Sigma$
is called \emph{star-shaped} if there exists $x \in U$ such that 
$U \subset U(x)$ and $\overline{xy} \subset U$ for all $y \in U$. 
Here, $\overline{xy}$ denotes the line segment between $x$ and $y$.

For the purposes of this section, a \emph{smooth tropical variety without points of higher sedentarity}
is a Hausdorff topological space $X$ together with a finite smooth tropical atlas
(c.f.\ \cite[Definition 6.1]{MR-TropicalGeometry, FR-DiagonalTropicalMatroid}), 
that is, 
\begin{itemize}
  \item a finite open cover $X = \bigcup_{i \in I} U_i$, 
	\item loopfree matroids $M_i$, $i \in I$, 
	\item maps $\phi_i \colon U_i \to \Sigma_{M_i} \subset \R^{N_i}$ which 
	      are homeomorphisms onto their images,
	\item the transition maps $\phi_i \circ \phi_j^{-1}$ are locally $\Z$-affine 
	      (that is, the restriction of an affine map $\R^{N_j} \to \R^{N_i}$ 
				whose differential is defined over $\Z$),
	\item for any $J \subset I$ there exists $i \in J$ such that 
	      $\phi_i(U_J)$ is star-shaped ($U_J = \bigcap_{j \in J} U_j$).
\end{itemize}
We assume that $\dim(X) =n$, that is, the dimension of the fans $\Sigma_{M_i}$ is $n$.
The star-shape condition is mostly for simplicity: It allows us to use simpler homology versions
and avoids technical comparisons between (more involved) homology versions. 
Moreover, such an atlas can be easily constructed for example 
under the common assumption that $X$ admits a (global) polyhedral structure, 
e.g.\ \cite[Definition 1.10]{MZ-TropicalEigenwaveIntermediate}
and \cite[Definition 2.2]{JRS-Lefschetz11Theorem}.

A \emph{tropical subcycle} of $X$ is a weighted closed subset whose restriction to a chart
is an open subset of a tropical subcycle in $\R^{N_i}$ (balanced polyhedral set). 
Let $Z_k(X)$ denote the group
of $k$-dimensional tropical subcycles. 
The intersection-theoretic side of the trace formula is based on the intersection
product
\[
  Z_k(X) \times Z_l(X) \to Z_{n-k-l}(X) 
\]
defined in \cite[section 6]{FR-DiagonalTropicalMatroid}. 

A map $\psi \colon X \to Y$ between two smooth tropical varieties 
is \emph{morphism} if in charts it is locally $\Z$-affine. 
If $\psi$ is proper, we have have a pushforward 
$\psi_* \colon Z_k(X) \to Z_k(Y)$ which is defined locally
using \cite[Construction 7.3]{AR-FirstStepsTropical}.
Let $\psi \colon X \to X$ be a proper tropical endomorphism. 
The tropical subcycles $\Delta$ and $\Gamma_\psi$ of $X \times X$ 
are defined as the pushforwards of $X$ along 
$x \mapsto (x,x)$ and $x \mapsto (x,\psi(x))$, respectively. 

\begin{definition} 
  Let $\psi \colon X \to X$ be a proper tropical endomorphism of a smooth tropical variety $X$ without
	points of higher sedentarity. The \emph{cycle of stable fixed points} of $\psi$ is 
	the zero-dimensional cycle $\Gamma_\psi \cdot \Delta$ (or rather, its projection to $X$). 
\end{definition}

\subsection{Tropical homology groups}

Tropical homology  $H_{p,q}(X)$  and cohomology $H^{p,q}(X)$ were defined in
\cite[Section 2.4]{MZ-TropicalEigenwaveIntermediate, IKMZ-TropicalHomology}
as homology with local coefficients using 
the so-called framing groups $F_p$.
The Borel-Moore and compact support variants $H^{\BM}_{p,q}(X)$ and $H_c^{p,q}(X)$
appear for example in \cite[Definition 2.7]{JRS-Lefschetz11Theorem}. 
We only consider real coefficients case here and hence
drop $\R$ from the notation throughout. 
Adapted to our needs, we recall a singular and a Čech approach to these groups.
We refer to the aforementioned papers for more details
and the comparison to other definitions.

\paragraph{Tropical framing groups $F_p$}

Let $\Sigma$ be a polyhedral fan in $\R^N$,
$S$ the standard $q$-simplex and 
and $s \colon S \to \Sigma$ a singular $q$-simplex. 
Following \cite[Definition 13]{IKMZ-TropicalHomology}, 
we define the \emph{$p$-th framing group} for $s$ by
\[
  F_p(s) := \sum_{\substack{\sigma \in \Sigma \\ \im(s) \subset \sigma}} \textstyle \bigwedge^p \sigma.
\]
Here, $\bigwedge^p \sigma$ denotes the subspace of $\bigwedge^p \R^N$ $p$-wedges of vectors
in $\sigma$. Note that $F_p(s) \neq \{0\}$ only if there exists a cone $\sigma$ containing the image of $s$
(in more standard notation, $F_p(s) = F_p(\tau)$ where $\tau$ is the smallest such cone). 
If $T \subset S$ is a face, then clearly $F_p(s) \subset F_p(s|_T)$. 
For any $x \in X$, we define $F_p(x)$ by regarding $x$ as a singular $0$-simplex.

Given a star-shaped open subset $U \subset \Sigma$ with centre $x \in U$, we define similarly
\[
  F_p(U) := \sum_{\substack{\sigma \in \Sigma \\ x \subset \sigma}} \textstyle \bigwedge^p \sigma.
\]
If $y$ is another centre for $U$, $x$ and $y$ lie in the relative interior of the same cone of
$\Sigma$ (otherwise, $x \notin U(y)$), which shows that the definition is independent of the choice of centre. 
More generally, if $V \subset U$ is another star-shaped open subset with centre $y$, we have
$F_p(V) \subset F_p(U)$.

%


\paragraph{Singular tropical homology}

Let $X$ be a smooth tropical variety. 
A \emph{small $q$-simplex} in $X$ is a singular $q$-simplex
$s \colon S \to X$ 
such that there exists a chart $(U, \psi, M)$ of $X$
with $\im(s) \subset U$.
Then $\psi \circ s$ is a singular $q$-simplex in $\Sigma_M$ and
we define $F_p(s) :=  F_p(\psi \circ s)$. 
If $(U', \psi', M')$ is another chart containing $\im(s)$, we
have a canonical isomorphism $F_p(\psi \circ s) \cong F_p(\psi' \circ s)$,
and we use this identification without explicit mention (so really, an element
in $F_p(s)$ is a set of vectors, one for each such chart). 

%

We denote by $C_q(X, F_p)$ the direct sum of the vector spaces 
$F_p(s)$ for all small $q$-simplices $s$
and set $C_*(X, F_p) = \sum_{q=0}^\infty C_k(X, F_p)$.
We use the notation $\sum_s \alpha_s s$, $\alpha_s \in F_p(s)$, for a vector in $C_q(X, F_p)$.
We define a boundary map $\partial \colon C_q(X, F_p) \to C_{q-1}(X, F_p)$
by setting
\[
  \partial(\alpha s) = \alpha \partial s, 
\]
where $\partial s$ denotes the usual singular boundary map. 
This formula is well-defined since for any face $T \subset S$,
$s|_T$ is again small and $F_p(s) \subset F_p(s|_T)$ 
as mentioned above. 
Clearly, $\partial^2 = 0$ on $C_*(X, F_p)$.
The \emph{(singular) tropical homology groups} of $X$ 
are given by
\[
  H_{p,q}(X) := H_{q}(C_*(X, F_p), \partial),
\]
while the dual complex gives us tropical cohomology groups $H^{p,q}(X)$.
Finally, extending the construction to the differential complex 
$C_*^{\BM}(X, F_p)$ of locally finite chains of small simplices, 
we obtain the Borel-Moore and compact support versions
$H^{\BM}_{p,q}(X)$ and $H_c^{p,q}(X)$.

%
%

If $\psi \colon X \to Y$ is a tropical morphism, we have induced pushforward maps
$\psi_* \colon H_{p,q}(X) \to H_{p,q}(Y)$. If $\psi$ is proper, we also have maps
$\psi_* \colon H^{\BM}_{p,q}(X) \to H^{\BM}_{p,q}(Y)$. 
They are defined in the usual way with the
help of the local multi\hyp{}differentials $d\psi_s : F_p(s) \to F_p(\psi \circ s)$
on the level of coefficients.
Analogously, there are pullback maps $H^{p,q}$ and $H_c^{p,q}$ ($\psi$ proper). 
The trace side of the trace formula consists of the \emph{graded} trace of the map
$\psi_* \colon H^{\BM}_{*,*}(X) \to H^{\BM}_{*,*}(X)$.
Here, graded means that the trace of the piece of degree $(p,q)$ is counted
with sign $(-1)^{p+q}$. 
As mentioned in \autoref{rem:otherhomology}, we could equally well use
the graded trace of $\psi^*$ on $H_c^{*,*}(X)$. However, the graded trace
on $H_{*,*}(X)$ (or, equivalently, on $H^{*,*}(X)$) is different in general. 

\paragraph{Čech tropical homology}

The star-shape condition imposed on our (finite) atlas $(U_i, \psi_i, M_i)_{i\in I}$ allows us 
to give an alternative description of tropical homology based on 
Čech (co)homology.
By assumption, for $J \subset I$ there exists $i \in J$ such that $\psi_i(U_J)$ is star-shaped
in $\Sigma_{M_i}$, and we set $F_p(U_J) := F_p(\psi_i(U_J))$. 
As above, we use the canonical identification $F_p(\psi_i(U_J)) \cong F_p(\psi_{i'}(U_J))$
for another admissible index $i'$ without explicit mention. By convention, $F_p(\emptyset) = \{0\}$. 
We set 
\[
  \CCC_q(X, F_p) = \bigoplus_{\substack{J \subset I \\ |J| = q+1}} F_p(U_J)
\]
and $\CCC_*(X, F_p) = \sum_{q=0}^{|I|-1} \CCC_q(X, F_p)$.
After choosing a total order on $I$, we have a well-defined 
Čech differential $d \colon \CCC_q(X, F_p) \to \CCC_{q-1}(X, F_p)$
using the inclusions $F_p(V) \subset F_p(U)$ for $V \subset U$ 
mentioned above.
The \emph{Čech tropical homology groups} of $X$ 
are given by
\[
  \HHH_{p,q}(X) := H_{q}(\CCC_*(X, F_p), d).
\]

\begin{proposition} \label{propComparisonSingularCech}
  There are canonical isomorphisms
	\[
	  H_{p,q}(X) \cong \HHH_{p,q}(X).
	\]
\end{proposition}

To prove the proposition, we start with the following observation.

\begin{lemma} \label{lemStarShaped}
  Let $U \subset \Sigma$ be a star-shaped open subset. Then
	\[
	  H_{p,q}(U) = \begin{cases} 
		               F_p(U) & \text{if } q=0, \\
								   0 & \text{otherwise.}
								 \end{cases}
	\]
\end{lemma}

\begin{proof}
  This is standard. 
  Let $x \in U$ be a centre of $U$. Given a singular $q$-simplex $s$,
	we may consider its cone $c(s)$ over $x$; note that $F_p(s) = F_p(c(s))$. 
	We obtain a chain homotopy
	$c \colon C_q(X, F_p) \to C_{q+1}(X, F_p)$ such that 
	$\partial \circ c + c \circ \partial = \id$ for $q > 0$ and
	$\partial \circ c + c \circ \partial = \id - \kappa$ for $q = 0$.
	Here, $\kappa \colon C_0(X, F_p) \to F_p(x) \subset C_0(X, F_p)$
	is given by $\kappa(\alpha [y]) = \alpha$. Since $F_p(U) = F_p(x)$, 
	this proves the claim.
\end{proof}

\begin{proof}[Proof of \autoref{propComparisonSingularCech}]
  We consider the double complex given by
	\[
	  K_{r,s}(X,F_p) = \bigoplus_{\substack{J \subset I \\ |J| = r+1}} C_s(U_J, F_p).
	\]
	Its differentials 
	\begin{align*} 
	  \partial \colon K_{r,s}(X,F_p) \to K_{r,s-1}(X,F_p), && d \colon K_{r,s}(X,F_p) \to K_{r-1,s}(X,F_p), 
	\end{align*}
	are given by the singular boundary map 
	$\partial$
	and the Čech differential $d$
	(using, of course, $C_s(V, F_p) \subset C_s(U, F_p)$ for $V \subset U$). 
	By \autoref{lemStarShaped}, the homology $H_{r,s}(K_{*,*}(X,F_p), \partial)$
	with respect to $\partial$ is zero for $s > 0$ and equal to $\CCC_r(X, F_p)$
	for $s=0$. 
	Similarly, it is straightforward to check that the homology $H_{r,s}(K_{*,*}(X,F_p), d)$
	with respect to $d$ is zero for $r > 0$
	and equal to $C_s(X, F_p)$ for $r=0$. Indeed, note that by definition 
	the canonical map $K_{0,s}(X,F_p) \to C_s(X, F_p)$ is surjective; moreover, its 
	kernel is the image of $d \colon K_{1,s}(X,F_p) \to K_{0,s}(X,F_p)$.
	It now follows by general double complex yoga (e.g.\ 
	\cite[Section II.8]{BT-DifferentialFormsAlgebraic})
	that both $H_{p,q}(X)$ and $\HHH_{p,q}(X)$ are canonically
	isomorphic to the total homology of the double complex. This finishes the proof. 
\end{proof}

\subsection{The global tropical Poincaré-Hopf theorem} \label{subs:GlobalPoincareHopf}

We start with a definition. 

\begin{definition} 
  The \emph{tropical Euler characteristic} of a smooth tropical variety $X$ without points of higher sedentarity 
	is 
	\[
	  \chi(X) := \sum_{p,q} (-1)^{p+q} \dim H_{p,q}(X).
	\]
	The \emph{local tropical Euler characteristic} at a point $x$ is $\chi(x) := \sum_{p=0}^n (-1)^p F_p(x)$.
\end{definition}

\begin{remark} 
  Of course, the definition also makes sense for more general tropical spaces
	(non-smooth or with points of higher sedentarity).
	As mentioned in \autoref{rem:otherhomology}, by ordinary and Poincaré duality
	\cite[Theorem 5.3]{JSS-SuperformsTropicalCohomology, JRS-Lefschetz11Theorem} we could equivalently use
	the variants $H^{\BM}_{p,q}(X)$, $H^{p,q}(X)$ or $H_c^{p,q}(X)$ in the definition. 
\end{remark}

Summarizing some well-known facts, the following lemma asserts that
\autoref{mainresult} agrees with the special case of 
\autoref{traceformula} for $\psi=\id$ and $X = \Sigma_M$.

\begin{lemma} \label{lem:specialcase}
  Let $M$ be a loopless matroid of rank $n+1$. Then
	\[
	  (-1)^n \beta(M) = \chi(\Sigma_M) = \chi(0).
	\]
\end{lemma}

\begin{proof}
  By \autoref{lemStarShaped}, 
	the only non-zero homology groups of $\Sigma_M$ are
	\begin{equation} 
		H_{p,0}(\Sigma_M) = F_p(\Sigma) = F_p(0).
	\end{equation}
	Hence the right hand side equality is clear. 
	Moreover, it was shown in \cite[Theorem 4]{Zha-OrlikSolomonAlgebra} that 
	\begin{equation} 
		\bigoplus_p F_p(\Sigma_M) \cong \bigoplus_p \OS^p(M),
	\end{equation}
	where the right hand side denotes the Orlik-Solomon algebra of $M$. 
	Finally, by e.g.\ \cite[Theorem 3.68]{OT-ArrangementsHyperplanes}
	\begin{equation} \label{eq:betaEuler} 
		(-1)^n \beta(M) = \sum_p (-1)^p \dim \OS^p(M),
	\end{equation}
	and hence the claim follows. Alternatively, the left hand side
	follows from \cite[Theorem 5.1]{Rau-TropicalLefschetzHopf}.
\end{proof}

We call $x \in X$ a \emph{vertex} of $X$ if 
there exists a chart $(U,\psi, M)$ such that $x \in U$, $\psi(x) = 0$,
and $M$ connected. The latter condition is equivalent to 
$\Sigma_M$ having trivial lineality space $\{0\}$, c.f.\ 
\cite[Lemma 2.3]{FR-DiagonalTropicalMatroid}.
Since out atlas is finite, 
the set of vertices of $X$, denoted by $\Ver(X)$, is finite. 
We can now restate \autoref{GlobalPoincareHopf} (equivalently, \autoref{traceformula} for
$\psi = \id$) in the following refined form.

\begin{theorem}[Global tropical Poincaré-Hopf theorem] \label{GlobalPoincareHopfRefined}
  The tropical Euler characteristic of a smooth tropical variety without points of higher sedentarity
	is equal to 
	\[
	  \chi(X) = \sum_{x \in \Ver(X)} \chi(x) = \deg \Delta^2.
	\]
\end{theorem}

\begin{proof}
  In the special case $X = \Sigma_M$, is follows from \cite[Lemma 5.1]{FR-DiagonalTropicalMatroid}
	that $\Delta^2 = 0$ when $\Sigma_M$ has non-trivial lineality space.
	Thus, back to the general case, \autoref{mainresult}, \autoref{lem:specialcase} and the locality
	of the tropical intersection product \cite[Section 6]{FR-DiagonalTropicalMatroid}
	imply $\chi(x) = 0$ for $x \notin \Ver(X)$ and 
	\[
	  \Delta^2 = \sum_{x \in \Ver(X)} \chi(x) \cdot x
	\]
	(under the natural isomorphism $\Delta \cong X$).
	This proves the right hand side equality.
	
%
%
	%
	%
%
	
	To prove the left hand side, 
	we use Čech tropical homology  $\HHH_{p,q}(X)$.
	Note that $\CCC_*(X,F_p)$ is finite-dimensional. 
	Hence, by \autoref{propComparisonSingularCech} 
	\[
	  \chi(X) := \sum_{p,q} (-1)^{p+q} \dim \HHH_{p,q}(X) 
		         = \sum_{\substack{J \subset I \\ J \neq \emptyset}} (-1)^{|J|-1} \sum_p (-1)^p F_p(U_J).
	\]
	Consider an individual $U = U_J$ with centre $x$. 
	By definition, 
	\[
	  \sum_p (-1)^p F_p(U) = \sum_p (-1)^p F_p(x) = \chi(x).
	\]
	So this value is non-zero if and only if $U$ contains a vertex 
	which then is its unique centre. 
	We conclude
	\[
	  \chi(X) = \sum_{x \in \Ver(X)} \chi(x) \sum_{\substack{J \subset I \\ x \in U_J}} (-1)^{|J|-1}
		        = \sum_{x \in \Ver(X)} \chi(x).
	\]
	This proves the left hand side equality.
\end{proof}

\begin{remark} 
  In the proof, we used that $\chi(x) \neq 0$ implies that $x$ is a vertex.
	We may thus rewrite the left hand side equality as $\chi(X) = \int_X \chi(x) dx$. 
	In fact, it is known that the beta invariant of a loopfree matroid $M$ is $0$ if and only 
	if $M$ is disconnected \cite[Theorem 7.3.2]{Whi-CombinatorialGeometries}.
	Hence, $\chi(x) \neq 0$ for any vertex $x$.

	More generally, the $0$-cycle 
	$\text{csm}_0(X) := \sum_x \chi(x) x$ is the 
	natural extension of CSM classes \cite{LRS-ChernSchwartzMacpherson}
	from matroid fans to smooth tropical varieties.
	Comparing to the classical case, 
  it is interesting that 
	the tropical Poincaré-Hopf theorem can be localised 
	at the vertices of $X$ (the stable fixed points of $\psi=\id$)
	and stated on the level of $0$-cycles as $\Delta^2 = \text{csm}_0(X)$.
	This is of course related to the
	fact that the tropical intersection product is defined 
	on the cycle level without the need to pass to rational equivalence. 
\end{remark}

\subsection{The tropical Weil trace formula}

In this subsection, we prove \autoref{introcurves}.
In honour of Weil's formula for algebraic curves, we call this special case the 
tropical Weil trace formula.

Throughout this section, $C$ denotes a connected smooth tropical curve. Additionally to the charts described
previously, we also allow points of higher sedentarity with local model $-\infty \in \T = \R \cup \{-\infty\}$.
We denote by $\Ver(C)$ the set of vertices including the subset $\Ver^\infty(C)$ of points of higher sedentarity. 
We call $C^\mm := C \setminus \Ver^\infty(C)$ the \emph{mobile part} of $C$. 
Note that $C$ is irreducible in the sense that the group of $1$-dimensional tropical subcycles is
$Z_1(C) = \Z C \cong \Z$. Regarding tropical homology, the extra convention for $x \in \Ver^\infty(C)$
is $F_0(x) = \R$ and $F_1(x) = \{0\}$.
In other words, $H_{0,q}$ (still) computes 
ordinary homology with constant coefficients $\R$ while $H_{1,q}$ computes essentially
the relative homology for the pair $\Ver^\infty(C) \subset C$. 

\begin{definition} 
  Let $f \colon C \to D$ be a proper tropical morphism of connected smooth tropical curves $C$ and $D$.
	The \emph{degree} of $\psi$ is the integer $\deg(\psi) \in \N \cup \{0\}$ such that
	$\psi_*(C) = \deg(\psi) \cdot D$.
	
	An \emph{open edge} of $C$ is a connected component $e$ of $C \setminus \Ver(C)$.
	Its \emph{local degree} $\deg_e(\psi)$ is the absolute value of the local 
	stretching factor $d\psi_x \in \Z$, $x \in e$.
	\end{definition}

Here, we refer to the pushforward of tropical cycles defined for example in 
\cite[Construction 7.3]{AR-FirstStepsTropical}. Since we are only interested in endomorphisms,
the following lemma focuses on this case. 
For $l \in \R_>$, consider the tropical curve $S^1_l = \R/l\Z$. 
It is a tropical elliptic curve without vertices whose loop has length $l$.

\begin{lemma} \label{endomprop}
  Let $f \colon C \to C$ be a proper tropical endomorphism of a connected smooth tropical curve $C$.
	Then the following holds. 
	\begin{enumerate}
		\item If $\deg(\psi) = 0$, $\psi$ is constant.
		\item If $\deg(\psi) > 0$, $\psi$ is surjective and $e \mapsto \psi(e)$ is a bijection
		      on the set of open edges. 
		\item If $X \not\cong S^1_l$, then $\deg_e(\psi) = \deg(\psi)$ for any open edge $e$. 
		\item If $\deg(\psi) = 1$, $\psi$ is an automorphism.
		\item If $\deg(\psi) \geq 2$, then either $C \cong S^1_l$, or $X^\mm \cong (-\infty, 0)$
		      or $X^\mm \cong \Sigma_M$, where $M$ is a loopfree matroid of rank $2$.
	\end{enumerate}	
\end{lemma}

\begin{proof}
  By definition of the pushforward $\psi_*(C)$ \cite[Construction 7.3]{AR-FirstStepsTropical},
	the degree can be computed at generic points of $C$ by counting preimages $x$ with (positive)
	weights $|d\psi_x|$. 
	Since $\psi(C) \subset C$ is a connected subgraph, this proves (a) and (since $\psi$ is proper) 
	the first part of (b). The set of open edges is finite and the image $\psi(e)$
	of an open edge is either a point or contained in an open edge (balancing condition). 
	Hence, surjectivity implies that
	in fact $\psi(e)$ is equal to an open edge and that this assignment is bijective. 
	
	An open edge is isometric to one the following four models: $(0,l), l \in \R_>$,
	$(0,+\infty)$, $\R$ or $S^1_l$. The restriction $\psi|_e : e \to \psi(e)$ is an affine, surjective map in
	the first three
	cases and therefore bijective. By the previous remarks, statement (c) follows. 
	
	If $\deg(\psi) = 1$, $\psi$ is invertible (over $\Z$) on $C \setminus \Ver(C)$ by (c)
	and it is clear that this can be extended to the vertices (note that it follows 
	from the previous statements that $\psi$ also induces a bijection on $\Ver^\infty(C)$). 
  This proves (d).
		
	For (e), note that an open edges isometric to $(0,l)$ must be mapped to an edge
	of type $(0, \deg_e(\psi) l)$. Hence, using (c), such edges cannot exist if $\deg(\psi) \geq 2$
	(take an edge of maximal length). Curves without such edges fall into one of the classes
	listed in (e). 
\end{proof}

\begin{remark} 
  We say $C$ is of finite type if every chart $(U,\phi)$ can be extended to a chart $(U',\phi')$,
	$U \subset U'$, such that $\overline{\phi(U)} \subset \phi'(U)$ (the closure is taken in $\R^N$)
	\cite[Definition 6.1.14]{MR-TropicalGeometry}. 
	For curves, this is equivalent to the requirement that 
	for any open edge $e$ isometric to $(0,l)$ or $(0,+\infty)$ the limit in $C$ for 
	$x \to 0$ exists (by symmetry, also for $x \to l$ in the first case). 
	For such curves, one can show that any tropical endomorphism $\psi \colon C \to C$ is
	either constant or proper (and hence surjective). Indeed, if $\psi$ is non-constant,
  one can show that $\psi(C)$ is open (using local irreducibility) and that $C^\mm \subset \psi(C)$
	(using finite type). Hence, again, $e \mapsto \psi(e)$ defines a bijection of open edges
	which can be restricted to those edges whose closure contains a point of higher sedentarity. 
	This implies surjectivity and properness. 
\end{remark}

Next, we discuss the local cases of computing $\Gamma_\psi \cdot \Delta$. 
Let $M$ be a loopfree matroid of rank $2$ and let $\psi : \Sigma_M \to \Sigma_M$ 
be a proper (i.e.\ non-constant) tropical endomorphism.
Without loss of generality, we may restrict the ambient space $\R^N$ to the span of $\Sigma_M$
and, equivalently, assume that all rank $1$ flats are singletons. 
Under this assumption, the permutation of rays of $\Sigma_M$ under $\psi$ induces a bijection $\psi' \colon E \to E$.
We denote by $\fix(\psi') = \#\Fix(\psi')$ the number of elements fixed by $\psi'$.

\begin{lemma} \label{curvelocalcase1}
  Let $M$ be a loopfree matroid of rank $2$ and let $\psi : \Sigma_M \to \Sigma_M$ 
	be a proper tropical endomorphism such that $\psi(0) = 0$. Then 
	\begin{equation} 
		\deg (\Gamma_\psi \cdot \Delta) = \deg(\psi) + 1 - \fix(\psi').	
	\end{equation}
\end{lemma}

\begin{proof}
  We set $d := \deg(\psi)$. 
	According to \autoref{prop:diagonalintersection}, we can compute 
	$\Gamma_\psi \cdot \Delta$ as $\gamma^*(g_1) \cdot \Sigma_M$, 
	where $g_1$ is the function from \autoref{eq:functiongidehom} and 
	$\gamma \colon x \mapsto (x, \psi(x))$.
	The image of the primitive generator $v_{\{i\}}$ under $\gamma$ 
	in terms of primitive generators for $\Sigma_{M \oplus_0 M}$ is
	\begin{equation*} 
		\gamma(v_{\{i\}}) = 
		            \begin{cases}
								 v_{(\{i\},\{\psi'(i)\})} + (d-1) v_{(\emptyset,\{\psi'(i)\})} & i \neq 0 \neq \psi'(i), \\
								 v_{(\{0\},\{0\})} + (d-1) v_{(E,\{0\})}                       & i = 0 = \psi'(i), \\
								 v_{(\{0\},E)} + d v_{(\emptyset,\{\psi'(i)\})}                & i = 0 \neq \psi'(i), \\
								 v_{(\{i\},\emptyset)} + d v_{(E,\{0\})}                       & i \neq 0 = \psi'(i).
							 \end{cases}
	\end{equation*}
	It follows that 
	\begin{equation*} 
		\gamma^*(g_1)(v_{\{i\}}) = 
		            \begin{cases}
								 -1 & i = \psi'(i) \neq 0, \\
								  0 & i \neq 0 \neq \psi'(i) \neq i, \\
								 +d & i = 0 = \psi'(i), \\
								 +1 & i = 0 \neq \psi'(i), \\
								 +d & i \neq 0 = \psi'(i).
							 \end{cases}
	\end{equation*}
	Since $\deg (\Gamma_\psi \cdot \Delta) = \sum_{i \in E} \gamma^*(g_1)(v_{\{i\}})$, the claim follows. 
\end{proof}

In the presence of points of higher sedentarity, we also have to compute
the contribution of such points to $\Gamma_\psi \cdot \Delta$. 
To do so, we use the extension of the intersection product to 
points of higher sedentarity for smooth tropical \emph{surfaces} 
\cite{Sha-TropicalSurfaces} (see also \cite{MR-TropicalGeometry}).
The only extra ingredient here is the assignment of an intersection multiplicity 
for two lines in $\T^2$ meeting in $(-\infty, -\infty)$. Given
the primitive generators $(a_1, a_2)$, $(b_1,b_2) \in \N^2$,
this multiplicity is set to be $\min\{a_1 b_2, a_2 b_1\}$
\cite[Definition 3.5.1]{MR-TropicalGeometry}.

\begin{lemma} \label{curvelocalcase2}
  Let $\psi : [-\infty, 0) \to [-\infty, 0)$ 
	be a proper tropical endomorphism. Then $\deg (\Gamma_\psi \cdot \Delta) = 1$.
\end{lemma}

\begin{proof}
  Any such $\psi$ is of the form $x \mapsto dx$ with $d = \deg(\psi)$. 
	The two cycles $\Delta$ and $\Gamma_\psi$ are two rays in 
	$[-\infty, 0)^2$ with primitive direction vectors $(1,1)$ and $(1,d)$. 
  Their intersection is equal to the point 
	$(-\infty,-\infty)$ with multiplicity
	$\min\{1 \cdot d, 1 \cdot 1\} = 1$. 
\end{proof}

\begin{remark} 
  Note that in the situation of the lemma the only non-zero 
	Borel-Moore homology group is
	$H^{\BM}_{1,1}(X) = \R$ and $\Tr(\psi_*, H^{\BM}_{1,1}(X)) = d$. 
	The discrepancy on the trace  side gets corrected if we extend $\psi$ to a map
	$\T \to \T$, since now there is an extra fixed point $0$ 
	with intersection multiplicity $d-1$ by \autoref{curvelocalcase1}.
	However, the special case $C^\mm \cong (-\infty,0)$ and $\deg(\psi) > 1$ must be excluded
	from the following trace formula. 
\end{remark}

To compute the trace side of our upcoming formula, 
it will be useful to use the cellular version of tropical homology 
\cite[Section 2.2]{MZ-TropicalEigenwaveIntermediate}:
A \emph{cell structure} $\ZZ =(Z_0, Z_1)$ for $C$ consists of
a non-empty finite subset $Z_0 \subset C$ containing $\Ver(C)$ and
$Z_1$ the set of connected components of $C \setminus Z_0$
such that the closure of every $e \in Z_1$ is homeomorphic to either $[0,1)$ or $[0,1]$.
We call the elements of $Z_0$ and $Z_1$ $0$-cells and $1$-cells, respectively. 
We set $C_q^\cell(\ZZ, F_p) = \bigoplus_{x \in Z_q} F_p(x)$. 
Choosing orientations for $1$-cells, we obtain a well-defined cellular boundary map 
$\partial_\ZZ \colon C_1^\cell(\ZZ, F_p) \to C_0^\cell(\ZZ, F_p)$
(when $p = 1$, boundary points $x$ of $e \in Z_1$ of higher sedentarity are discarded 
since $F_1(x) = \{0\}$). 
In fact, $\partial_\ZZ$ can be easily described as follows. 
The relative homology groups $H^\BM_{p,q}(\overline{e}, \partial e; F_p)$
are zero for $q \neq 1$ and canonically isomorphic (using the orientations) to $F_p(e)$ for $q=1$. 
Hence, $H^\BM_{1}(C, Z_0; F_p) \cong C_1^\cell(\ZZ, F_p)$ and 
one checks easily that the long exact sequence associated to the pair $Z_0 \subset C$
boils down to the exact sequence 
\[
  0 \rightarrow H^\BM_{p,1}(X) \rightarrow 
	  C_1^\cell(\ZZ, F_p) \stackrel{\partial_\ZZ}{\rightarrow} C_0^\cell(\ZZ, F_p) 
		\rightarrow H^\BM_{p,0}(X) \rightarrow 0.
\]
In particular, cellular homology agrees with singular Borel-Moore homology. 

Given a proper tropical  endomorphism $\psi \colon C \to C$, assume that
there exists a cell structure $\ZZ$ such that $\psi(c) \in \ZZ$ for $c \in \ZZ$. 
Then the pushforward on homology is induced by the map 
$\psi_* \colon C_q^\cell(\ZZ, F_p) \to C_q^\cell(\ZZ, F_p)$
given for $c \in Z_q$, $\alpha \in F_p(c)$, $x \in e$ by
\[
  \psi_*(\alpha c) = d \psi_x (\alpha) \psi(c)
\]
when $\psi(c) \in \ZZ_q$ and $\psi_*(\alpha c) = 0$ otherwise. 
The classical Hopf trace lemma \cite[§9, Theorem 2.1]{GD-FixedPointTheory}
states that the graded trace on homology can then be 
computed as the graded trace on the chain level,
\begin{equation} \label{eq:HopfTraceLemma} 
  \sum_{p,q} (-1)^{p+q} \Tr(\psi_*, H^{\BM}_{p,q}(C)) = \sum_{p,q} (-1)^{p+q} \Tr(\psi_*, C_q^\cell(\ZZ,F_p)).
\end{equation}
We are now ready to prove the main theorem for curves.

\begin{theorem}[Tropical Weil trace formula] 
  Let $\psi \colon C \to C$ be a proper tropical endomorphism of a connected smooth tropical curve $C$ 
	such that $C^\mm \not\cong (-\infty,0)$ or $\deg(\psi) = 1$. 
	Then we have
	\begin{equation} \label{eq:traceformulacurves} 
		\deg (\Gamma_\psi \cdot \Delta) = \sum_{p,q} (-1)^{p+q} \Tr(\psi_*, H^{\BM}_{p,q}(C)).
	\end{equation}
\end{theorem}

\begin{proof}
  Let us first deal with a few special cases. 
	If $\deg(\psi) = 0$, then by part (a) of \autoref{endomprop},  
	$\psi \equiv c$ is constant and 
	$\deg(\Gamma_\psi \cdot \Delta) = \deg(C \cdot \{c\}) = 1$. 
	So, \autoref{eq:traceformulacurves} holds true after replacing $H^{\BM}_{p,q}$ with $H_{p,q}$.
	Since $\psi$ is proper, $C$ is compact which implies $H^{\BM}_{p,q} = H_{p,q}$.
	
	The case $X \cong S^1_l$ is covered by \autoref{traceformulatori},
	so we exclude this case here. 
		
	Let us now assume $\deg(\psi) \geq 2$. By \autoref{endomprop}  (e) 
	and the exclusions made so far, this implies $C^\mm = \Sigma_M$
	where $M$ is a loopfree matroid of rank $2$. 
	Note that we can assume $\psi(0) = 0$, since even in the case $C^\mm = \R$
	there exists a fixed point in $C^\mm$, given that $\deg(\psi) \geq 2$.
	We use the same assumptions
	and notation as in \autoref{curvelocalcase1}.
	Additionally, we denote by $\fix(\psi^\infty)$ the number of 
	points in $\Ver^\infty(C)$ fixed by $\psi$. 
	\autoref{curvelocalcase1} and \autoref{curvelocalcase2} imply
	that $\deg(\Gamma_\psi \cdot \Delta) = \deg(\psi) + 1 - \fix(\psi') + \fix(\psi^\infty)$.
	
	To compute the trace side, 
	we use cellular homology for the cell structure $\ZZ$ given by 
	$Z_0 = \Ver(C)$.
	We see directly that $H^{\BM}_{1,0}(C) = 0$, $H^{\BM}_{1,1}(C) = \R C$
	and $\Tr(\psi_*, H^{\BM}_{1,1}(C)) = \deg(\psi)$. For $p=0$, 
	we use the Hopf trace lemma \autoref{eq:HopfTraceLemma}.
	We check easily that $C^{\cell}_{1}(\ZZ, F_0) = \R^{Z_1}$ and 
	$\Tr(\psi_*, C^{\cell}_{1}(\ZZ, F_0)) = \fix(\psi')$
	while $C^{\cell}_{0}(\ZZ, F_0) = \R^{Z_0}$ and 
	$\Tr(\psi_*, C^{\cell}_{0}(\ZZ, F_0)) = 1 + \fix(\psi^\infty)$, as required.
	
	We are left with the case $\deg(\psi) = 1$, so $\psi$ is an automorphism by
	\autoref{endomprop} (d). 
	By subdividing the open edges isometric to $(0,l)$ or $\R$ which get flipped
	by $\psi$, we can construct a cell structure  $\ZZ$ such that
	$\psi(c) \in \ZZ$ for $c \in \ZZ$ and $\psi|_e = \id$ whenever $\psi(e) = e$. 
	In particular, $\Fix(\psi)$ is a union of cells in $\ZZ$.
	Let us denote by $\psi_0$ and $\psi_1$ the induced bijections on $0$-cells and $1$-cells,
	respectively. Moreover, for any $x \in \Fix(\psi_0)$ we denote by $\psi_x$ 
	the permutation of the $1$-cells containing $x$. 
	Using \autoref{curvelocalcase1} and \autoref{curvelocalcase2} again, we get
	\begin{equation*} 
	  \deg(\Gamma_\psi \cdot \Delta) = \fix(\psi^\infty) + \sum_{\substack{x \in \Fix(\psi_0)\\ x \notin \Ver^\infty(C)}} 2 - \fix(\psi_x) 
		                               = \sum_{x \in \Fix(\psi_0)} 2 - \fix(\psi_x).
	\end{equation*}
	Setting $\Tr_{p,q} := \Tr(\psi_*, C^{\cell}_{q}(\ZZ, F_p))$, we have
	\begin{align} 
		\nonumber  \Tr_{0,0} &= \fix(\psi_0),     &   \Tr_{1,0} &= \sum_{x \in \Fix(\psi_0)} \fix(\psi'_x) - 1, \\
		\nonumber  \Tr_{0,1} &= \fix(\psi_1),     &   \Tr_{1,1} &= \fix(\psi_1).
	\end{align} 
  Indeed, the cases $p=0$ are obvious. In the cases $p = 1$, we use $\deg(\psi) = 1$. 
	Moreover, considering $\Tr_{1,0}$, the contribution of each $x \in \Fix(\psi_0) \setminus \Ver^\infty(C)$ is equal to $\Tr(d\psi_x, F_1(x))$.
	With the assumption from \autoref{curvelocalcase1}, we can resolve $F_1(x)$ by
	\[
	  0 \to \R\mathbf{1} \hookrightarrow \R^E \to F_1(x) \to 0.
	\]
	It follows that $\Tr(d\psi_x, F_1(x)) = \fix(\psi_x) - 1$. Finally, a point $x \in \Fix(\psi_0) \cap \Ver^\infty(C)$
	contributes zero to $\Tr_{1,0}$ since $\F_1(x) = 0$. 
	Finally, using \autoref{eq:HopfTraceLemma} again, we get  
	\[
	  \sum_{p,q} (-1)^{p+q} \Tr(\psi_*, H^{\BM}_{p,q}(C)) = \fix(\psi_0) - \sum_{x \in \Fix(\psi_0)} (\fix(\psi'_x) - 1) = \sum_{x \in \Fix(\psi_0)} 2 - \fix(\psi'_x).
	\]
	This proves the claim.
\end{proof}

\begin{remark} \label{workswithusualhom}
  In fact, in the case $\deg(\psi) \leq 1$ (or rather, $\psi$ constant or automorphism, 
	even without the properness assumption), our proof works equally well with 
	usual homology $H_{p,q}(C)$. Indeed, the only change required is that now 
	$\Tr_{0,1} = \Tr_{1,1} = \fix(\psi_1^b)$ only counts \emph{compact} fixed edges.
\end{remark}

\begin{example} 
  Let $X$ be a tropical curve whose underlying graph is the $\Theta$-graph $G$
	and with vertices $v_1$ and $v_2$ (see \autoref{exampleTheta}). 
	Its homology groups can be easily calculated as $H_{0,0}(X) =\Z$, 
	$H_{1,0}(X) = F_1(v_1) \cong \Z^2$, $H_{0,1}(X) = H_1(G) \cong \Z^2$
	and $H_{1,1}(X) = \Z\cdot X = \Z$.
	
	Let $\psi_1 \colon X \to X$ be the automorphisms which exchanges $v_1$ and $v_2$
	and flips every edge. Then $\Gamma_{\psi_1}\cdot \Delta$ consists of the midpoints
	of the three edges, and each of them occurs with intersection multiplicity $2$ by
	\autoref{curvelocalcase1}. On the trace side, $\psi_1$ induces $\id$ for
	$(p,q) = (0,0)$ and $(1,1)$ and $-\id$ for $(1,0)$ and $(0,1)$. 
	So the trace side gives $1 -(-2) -(-2) +1 = 6$ as well.
	
	Assume now that two of the edges of $X$ have the same length. Then there exists
	an automorphism $\psi_2 \colon X \to X$ which exchanges the two edges (but keeps
	the vertices and the third edge fixed). The set-theoretic fixed point locus of
	$\psi_2$ consists of the third edge, but only the vertices are stable fixed points. 
	Each vertex has intersection multiplicity $1$ in $\Gamma_{\psi_2}\cdot \Delta$.
	The pushforward is still identity for $(0,0)$ and $(1,1)$. For $(1,0)$ and $(0,1)$, however,
	we can choose bases such that $(\psi_2)_*$ permutes the basis elements. 
	Hence these traces are zero and we get $1 + 0+0+1 =2$.	
\end{example}

	\begin{figure}[!h]
		\centering
		\input{pic/exampleTheta.TpX}
		\caption{The two endomorphisms $\psi_1$ and $\psi_2$ and their stable fixed points with multiplicities}
		\label{exampleTheta}
	\end{figure}%

\subsection{The trace formula for tropical tori}

Let $\Lambda \subset \R^n$ be a lattice in $\R^n$ (i.e.\ a discrete free abelian subgroup).
The quotient $X = \R^n/\Lambda$ is called a \emph{tropical torus}. Note that the \enquote{tropical} (here, integral-affine)
structure on $X$ is induced by the lattice $\Z^n \subset \R^n$.
For more information on tropical tori, we refer to \cite[Section 5]{MZ-TropicalCurvesTheir}. 
We denote by $\Mat(\Lambda)$ the $n \times n$ matrices with real entries such that $A \Lambda \subset \Lambda$.

\begin{lemma} 
  Let $X = \R^n/\Lambda$ be a tropical torus and $\psi \colon X \to X$ a tropical endomorphism. Then $\psi$
	is of the form 
	\[
	  x \mapsto Ax + v \mod \Lambda
	\]
	for some matrix $A \in \Mat(\Z^n) \cap \Mat(\Lambda)$ and $v \in \R^n$. Moreover,
	$A$ is uniquely determined by $\psi$. 
\end{lemma}

\begin{proof}
  For any $x \in X$, we can canonically identify $T_x X = \R^n$ with 
	lattice of integer tangent vectors $T_x^\Z X = \Z^n$.
	By definition of tropical morphisms (see e.g.\ \cite[Section 7.3]{MR-TropicalGeometry}),
	the differential map $x \mapsto d\psi_x \colon T_x X \to T_x X$ is locally constant
	with values in $\Mat(\Z^n)$. Hence, it is globally constant, and setting 
	$A = d\psi_x$ and $v$ such that $[v] = \psi(0)$, we see that $\psi$ the required form. 
	The compatibility with $\Lambda$ implies $A \in \Mat(\Lambda)$. 
\end{proof}

\begin{theorem} \label{traceformulatori}
  Let $X = \R^n/\Lambda$ be a tropical torus and $\psi \colon X \to X$ a tropical endomorphism
	with differential $d\psi = A \in \Mat(\Z^n) \cap \Mat(\Lambda)$.
	Then 
	\begin{equation} \label{eq:traceformulatori} 
		\deg (\Gamma_\psi \cdot \Delta) = \det(\id - A)^2 = \sum_{p,q} (-1)^{p+q} \Tr(\psi_*, H_{p,q}(X)) .
	\end{equation}
\end{theorem}

\begin{proof}
  Let us first do the linear algebra behind the statement.
  Let 
	\[
	  \chi(t) = \det(\id - t A)
	\]
	be the characteristic polynomial of $A$ (for $s = 1/t$). 
	By standard expansion of determinants, the coefficient of $(-t)^k$ in $\chi$
	is equal to the sum of minors $\det(A_I)$ of size $k$.
	Here, $A_I$ denotes the diagonal submatrix of $A$ with rows and columns given by
	$I \subset \{1,\dots,n\}$. 
	On the other hand, consider the map $A^{\wedge k} \colon \bigwedge^k \R^n \to \bigwedge^k \R^n$ induced
	by $A$. Clearly, its trace is also equal to the sum of minors $\det(A_I)$ of size $k$. 
	(A fancier way of saying the same thing is that $\chi(t)$ is equal to the \emph{graded} trace 
	(i.e.\ graded pieces are counted with alternating signs) of the map
	$tA^{\wedge *} \colon \bigwedge^* \R^n \to \bigwedge^* \R^n$ induced
	by $tA$ on the exterior algebra.)
	
	Now, since the framing groups $F_p(x)$ are constant on tropical tori, 
	the $(p,q)$ groups are equal to 
	\[
	  H_{p,q}(X) = F_p([0]) \otimes H_q(X,\R) = \bigwedge^p \R^n \otimes \bigwedge^q \R^n.
	\]
	Moreover, under this identification, $\psi_*$ is equal to $A^{\wedge *} \otimes A^{\wedge *}$.
	By the computation from above the graded trace of $\psi_*$ is hence equal to
	\[
	  \Tr(\psi_*) = \Tr(A^{\wedge *})^2 = \det(\id - t A)^2,
	\]
	which proves the right hand side of \autoref{eq:traceformulatori}. 
	
	It remains to check the left hand side. 
	Consider the subspaces $D = \{(x,x)\}$ and $G = \{x, Ax\}$ of $\R^{2n}$. 
	Clearly, $\Delta$ and $\Gamma_\psi$ are translations of the projections 
	of $D$ and $G$, respectively, modulo $\Lambda \times \Lambda$. 
	Hence, for any fixed point of $\psi$ we can compute its contribution
	to $\Gamma_\psi \cdot \Delta$ as the tropical intersection multiplicity
	of $G$ and $D$. This multiplicity can be computed by combining lattice bases
	of $G$ and $D$ in a matrix and taking the absolute value of its determinant,
	see \cite[Lemma 5.1]{MR-TropicalDescendantGromov}.
	In our case, we get the absolute value of 
  \[
	  \det\left(
			\begin{array}{c|c}
				\id & \id \\ \hline
				A & \id
			\end{array}\right)
	  = \det\left(
			\begin{array}{c|c}
				\id & 0 \\ \hline
				A & \id - A
			\end{array}\right)
		= \det(\id - A).
	\]
	If $\det(\id - A) = 0$, this implies $\Gamma_\psi \cdot \Delta = 0$ since 
	all intersection multiplicities are zero. 
	Now assume $\det(\id - A) \neq 0$. In this case, the equation
	$(\id - A)x = v$ has exactly $|\det(\id - A)|$ many solutions modulo $\Lambda$ 
	(since $|\det(\id - A)|$ also computes the degree of the map 
	$\id - A \colon \R^n/\Lambda \to \R^n/\Lambda$).
	Hence $\psi$ has exactly $|\det(\id - A)|$ fixed points, and each contributes
	$|\det(\id - A)|$ to $\deg(\Gamma_\psi \cdot \Delta)$, which proves
	the left half of \autoref{eq:traceformulatori}.	
\end{proof}

\begin{remark} 
  Note that the tropical formula is essentially a product 
	of the \emph{classical} trace formula for $\psi$ (considered as map between manifolds)
	with the tropical formula for the map $A \colon \R^n \to \R^n$.
	Indeed, from the classical point of view, each of the $|\det(\id - A)|$ fixed points
	intersect transversally, and the homology groups to consider are 
	$H_*(X,\R) = \bigwedge^* \R^n$. 
	The first part of the proof is then a proof of the classical version
	(when using the correct signs).
	On the other hand, the tropical homology groups for $\R^n$ are
	$H_{*,*}(\R^n) = F_*(\R^n) = \bigwedge^* \R^n$ as well. 
	Again, the correct pieces of the previous proof also prove
	the trace formula for $A \colon \R^n \to \R^n$.
	Finally, combining the two parts gives the \enquote{squared} version
	of \autoref{traceformulatori}.
\end{remark}

\printbibliography


\subsection*{Contact}

    Johannes Rau \\ Departamento de Matemáticas \\ Universidad de los Andes \\
		KR 1 No 18 A-10, BL H \\ Bogotá, Colombia \\ \href{mailto:j.rau@uniandes.edu.co}{j.rau AT uniandes.edu.co}

\end {document}